# A study on the transmission dynamics of COVID-19 considering the impact of asymptomatic infection

Zonghao Zhang, Xiaotong Huang, Kedeng Cheng, Chuanqing Xu*, Songbai Guo, Xiaojing Wang
School of Science
Beijing University of Civil Engineering and Architecture
Beijing 102699, China

**Abstract:** The COVID-19 epidemic has been spreading around the world for nearly three years, and asymptomatic infections have exacerbated the spread of the epidemic. To evaluate the role of asymptomatic infections in the spread of the epidemic, we develop mathematical models to assess the proportion of asymptomatic infections caused by different strains of the main covid-19 variants. The analysis shows that when the control reproduction number is less than 1, the disease-free equilibrium point of the model is globally asymptotically stable; and when the control reproduction number is greater than 1, the endemic equilibrium point exists and is unique, and is locally asymptotically stable. We fit the epidemic data in the four time periods corresponding to the selected 614G, Alpha, Delta and Omicron variants. The fitting results show that, from the comparison of the four time periods, the proportion of asymptomatic persons among the infected persons gradually increased. We also predict the peak time and peak value for the four time periods, and the results indicate that the transmission speed and transmission intensity of the variant strains increased to some extent. Finally, we discuss the impact of the detection ratio of symptomatic infections on the spread of the epidemic. The results show that with the increase of the detection ratio, the cumulative number of cases has dropped significantly, but the decline in the proportion of asymptomatic infections is not obvious. Therefore, in view of the hidden transmission of asymptomatic infections, the cooperation between various epidemic prevention and control policies is required to effectively curb the spread of the epidemic.

1. Introduction

Since its outbreak at the end of 2019, COVID-19 has had a severe impact on countries around the world. How to effectively control the spread of the epidemic and restore normal production and life is still an issue that governments need to consider. To control the spread of COVID-19, many measures have been taken, including non-pharmacological interventions, pharmacological interventions, vaccinations, and more. However, due to the continuous emergence of variants, the epidemic prevention measures did not achieve the expected effect.

---

* Corresponding Author: xuchuanqing@bucea.edu.cn

There are a large number of contagious asymptomatic infections during the spread of COVID-19 [1,2,3,4]. This feature determines the high concealment of the spread of COVID-19, and leads to an increase in the difficulty of epidemic prevention and control. There have been a large number of research results on asymptomatic infections of new coronary pneumonia, including age distribution, average proportion, and transmission intensity [5,6,7]. In addition, some researchers have studied the spread of COVID-19 by establishing dynamic models including asymptomatic infections. Ruan et al. established a time-varying COVID-19 transmission compartment model including asymptomatic infected persons, simulated and reviewed the development process of the Wuhan epidemic, and obtained that the asymptomatic proportion of infected persons was about 20% [8]. Rahul Subramanian et al. established a COVID-19 transmission model including asymptomatic infections, quantified asymptomatic infections in New York City, and obtained asymptomatic infections accounted for about 60% [9]. Mohamed Amouch et al. proposed a new epidemiological mathematical model of the spread of COVID-19 disease, fitted the outbreak in Monaco, and obtained that the proportion of asymptomatic patients was 30% [10]. However, these studies are all conducted for a specific variant, and cannot effectively reflect the changes caused by variant iterations.

In the more than two years since the first appearance of COVID-19, many countries and regions around the world have experienced repeated outbreaks. Taking England as an example, although multiple rounds of lockdown measures have been adopted to control the spread of the epidemic, such epidemic prevention measures have not been fully effective. When the epidemic prevention measures were gradually lifted, the epidemic rebounded again. We believe that asymptomatic infections have played a very important role in the rebound of the epidemic. The lockdown measures have effectively reduced the number of infected people to a certain extent, but because large-scale screening tests were not adopted, there were still a certain number of asymptomatic infections in the population. After the epidemic prevention measures are lifted, these undetected asymptomatic infections will cause the next round of outbreaks. If the number of asymptomatic infections can be effectively estimated in the early stage of the next outbreak, and certain epidemic prevention and control measures are taken, the spread of the epidemic can be delayed to a certain extent. We establish an improved SEIAR infectious disease dynamics model to assess the role of asymptomatic infections in the early stages of epidemic transmission.

The remainder of this article is organized as follows. In Section 2, we present the COVID-19 infectious disease model and analyze the stability of the equilibrium point. In the third section, we use the established model to fit the actual epidemic data, and conduct a comparative analysis of the relevant kinetic parameters. In the last section, we conclude and discuss.

## 2. Model establishment and analysis

### 2.1 Model establishment

There are two typical characteristics of the spread of COVID-19. First, there are a large number of asymptomatic infected people with new coronary pneumonia. The second is that COVID-19 patients have a longer exposure period and have differences in infectivity. Based on our knowledge of COVID-19, we have the following assumptions in our model:

(A1) We divide the exposure period into two parts: the early stage and the later stage, in which the early stage is not infectious and the later stage is infectious;

(A2) For symptomatic infections, we assume that symptomatic infections can be screened as long as they are tested, and that the detected symptomatic infections are completely isolated and no longer contagious;

(A3) Asymptomatic infections were assumed not to be tested, and deaths from illness in asymptomatic infections were not considered.

Based on the above assumptions, we establish a COVID-19 transmission compartment model including asymptomatic infected persons. The model divides the general population into susceptible ($S$), pre-exposure patients ($E_1$), late-exposure patients ($E_2$), detected symptomatic infections ($I_1$), undetected symptomatic infections ($I_2$), asymptomatic infections ($A$) and recovered ($R$). Its dynamic flow chart is shown in Figure 1.

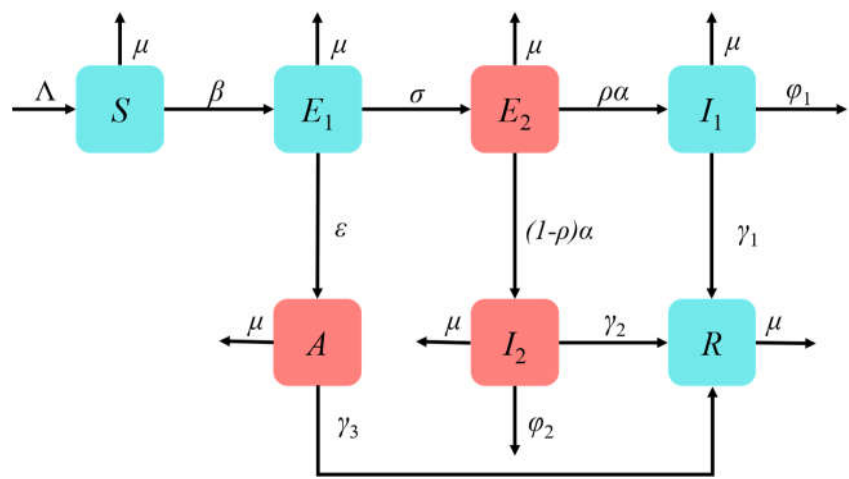

**Fig.1 Flow chart of the COVID-19 transmission dynamics model containing seven compartments.** The compartments represented by the red box is infectious, and the compartments represented by the blue box is not infectious.

The corresponding propagation dynamics equation is constructed as follows:

$$\begin{cases} \dfrac{dS}{dt} = \Lambda - \beta S(E_2 + I_2 + \omega A) - \mu S, \\ \dfrac{dE_1}{dt} = \beta S(E_2 + I_2 + \omega A) - (\sigma + \varepsilon + \mu)E_1, \\ \dfrac{dE_2}{dt} = \sigma E_1 - (\alpha + \mu)E_2, \\ \dfrac{dI_1}{dt} = \rho \alpha E_2 - (\gamma_1 + \varphi_1 + \mu)I_1, \\ \dfrac{dI_2}{dt} = (1-\rho)\alpha E_2 - (\gamma_2 + \varphi_2 + \mu)I_2, \\ \dfrac{dA}{dt} = \varepsilon E_1 - (\gamma_3 + \mu)A, \\ \dfrac{dR}{dt} = \gamma_1 I_1 + \gamma_2 I_2 + \gamma_3 A - \mu R. \end{cases} \quad (1)$$

The parameters in system (2.1) are explained as follows. $\Lambda$ is the daily number of births and $\mu$ is the daily natural mortality rate. The parameter $\beta$ represents the basal transmission rate, and $\sigma$ represents the conversion rate from patients in the pre-exposure period to patients in the late-exposure period. $\varepsilon$ represents the conversion rate from patients in the pre-exposure period to asymptomatic infected persons, and $\alpha$ represents the conversion rate from patients in the late exposure period to infected persons. $\omega$ is the weight of the transmission intensity of the asymptomatic infection relative to the symptomatic infection, and $\rho$ is the rate of testing among the symptomatic infection. $\gamma_1, \gamma_2, \gamma_3$ is the recovery rate of detected symptomatic infections, undetected symptomatic infections, and asymptomatic infections, respectively. $\varphi_1, \varphi_2$ is the mortality rate of detected symptomatic infections and undetected symptomatic infections, respectively.

2.2 Calculation of disease-free equilibrium point and controlled reproduction number

Obviously, there is always a disease-free equilibrium point $P_0 = (S_0, 0, 0, 0, 0, 0, 0)$ in system (1), where $S_0 = \Lambda/\mu$. Then we use the next generation matrix method to calculate the control reproduction number of system (1). The Jacobian matrices $F$ and $V$ at the disease-free equilibrium point are obtained from system (1) as

$$F = \begin{pmatrix} 0 & \beta S_0 & 0 & \beta S_0 & \omega\beta S_0 \\ 0 & 0 & 0 & 0 & 0 \\ 0 & 0 & 0 & 0 & 0 \\ 0 & 0 & 0 & 0 & 0 \\ 0 & 0 & 0 & 0 & 0 \end{pmatrix},$$

$$V = \begin{pmatrix} \sigma+\varepsilon+\mu & 0 & 0 & 0 & 0 \\ -\sigma & \alpha+\mu & 0 & 0 & 0 \\ 0 & -\rho\alpha & \gamma_1+\varphi_1+\mu & 0 & 0 \\ 0 & -(1-\rho)\alpha & 0 & \gamma_2+\varphi_2+\mu & 0 \\ -\varepsilon & 0 & 0 & 0 & \gamma_3+\mu \end{pmatrix}.$$

The control reproduction number is the spectral radius of the $FV^{-1}$, therefore we get

$$\mathcal{R}_C = \frac{\beta S_0}{\sigma+\varepsilon+\mu}\left(\frac{\sigma}{\alpha+\mu} + \frac{\sigma(1-\rho)\alpha}{(\alpha+\mu)(\gamma_2+\varphi_2+\mu)} + \frac{\varepsilon\omega}{\gamma_3+\mu}\right).$$

2.3 Existence of the endemic equilibrium point

Theorem 2.3. When $\mathcal{R}_C > 1$, the system (1) has a unique positive equilibrium point $P^* = (S^*, E_1^*, E_2^*, I_1^*, I_2^*, A^*, R^*)$.

Proof. Using the equilibrium equation of system (1) at the endemic equilibrium point, we can get

$$S^* = \frac{\Lambda}{\beta E_1^*\left(\frac{\sigma}{\alpha+\mu} + \frac{\sigma(1-\rho)\alpha}{(\alpha+\mu)(\gamma_2+\varphi_2+\mu)} + \frac{\varepsilon\omega}{\gamma_3+\mu}\right) + \mu},$$

$$E_2^* = \frac{\sigma}{\alpha+\mu}E_1^*,\ I_1^* = \frac{\sigma\rho\alpha}{(\alpha+\mu)(\gamma_1+\varphi_1+\mu)}E_1^*, \qquad (2)$$

$$I_2^* = \frac{\sigma(1-\rho)\alpha}{(\alpha+\mu)(\gamma_2+\varphi_2+\mu)}E_1^*,\ A^* = \frac{\varepsilon}{\gamma_3+\mu}E_1^*,$$

$$R^* = \frac{E_1^*}{\mu}\left[\frac{\sigma\rho\alpha\gamma_1}{(\alpha+\mu)(\gamma_1+\varphi_1+\mu)} + \frac{\sigma(1-\rho)\alpha\gamma_2}{(\alpha+\mu)(\gamma_2+\varphi_2+\mu)} + \frac{\varepsilon\gamma_3}{\gamma_3+\mu}\right].$$

Furthermore, according to the equilibrium equation, we obtain

$$\Lambda - \mu S^* - (\sigma+\varepsilon+\mu)E_1^* = 0.$$

Simplify to get

$$E_1^*\left[\beta(\sigma+\varepsilon+\mu)\mathcal{R}_C E_1^* + (\mu\beta S_0 - \Lambda\beta\mathcal{R}_C)\right] = 0.$$

Therefore, we can get

$$E_1^* = 0,\ \text{or}\ E_1^* = -\frac{\mu\beta S_0 - \Lambda\beta\mathcal{R}_C}{\beta(\sigma+\varepsilon+\mu)\mathcal{R}_C},$$

where

$$\mu\beta S_0 - \Lambda\beta\mathcal{R}_C = \Lambda\beta\left(\frac{\mu S_0}{\Lambda} - \mathcal{R}_C\right) = \Lambda\beta(1 - \mathcal{R}_C).$$

Therefore, $E_1^*$ has a unique positive root when $\mathcal{R}_C > 1$. According to equation (2), we can know that there is a unique positive equilibrium point $P^*$ of system (1) at this time.

2.4 Stability at equilibrium point

Theorem 2.4.1. When $\mathcal{R}_C < 1$, the disease-free equilibrium point $P_0$ of system (1) is locally asymptotically stable.

Proof. Let the Jacobian matrix of system (1) at the disease-free equilibrium point $P_0$ be $M_1$.

When $\mathcal{R}_C < 1$, assume that there is $\lambda_0$ satisfying $\text{Re}(\lambda_0) \geq 0$, such that $\det(\lambda_0 E - M_1) = 0$.
Expanding the characteristic equation, we can obtain

$$\begin{aligned}&(\lambda_0+\mu)^2(\lambda_0+\alpha+\mu)(\lambda_0+\gamma_1+\varphi_1+\mu)(\lambda_0+\gamma_2+\varphi_2+\mu)(\lambda_0+\gamma_3+\mu)\\&\{\lambda_0+(\sigma+\varepsilon+\mu)[1-F(\lambda_0)]\}=0,\end{aligned} \quad (3)$$

where

$$F(\lambda_0) = \frac{\beta S_0}{\sigma+\varepsilon+\mu}\left[\frac{\sigma}{\lambda_0+\alpha+\mu}+\frac{\sigma(1-\rho)\alpha}{(\lambda_0+\alpha+\mu)(\lambda_0+\gamma_2+\varphi_2+\mu)}+\frac{\varepsilon\omega}{\lambda_0+\gamma_3+\mu}\right].$$

The characteristic equation (3) is equivalent to

$$\lambda_0 + (\sigma+\varepsilon+\mu)[1-F(\lambda_0)] = 0.$$

The necessary condition for the above formula to be established is

$$\text{Re}[F(\lambda_0)] \geq 1. \quad (4)$$

Furthermore, since we assume $\text{Re}(\lambda_0) \geq 0$, we can obtain

$$\begin{aligned}|F(\lambda_0)| &= \left|\frac{\beta S_0}{\sigma+\varepsilon+\mu}\left[\frac{\sigma}{\lambda_0+\alpha+\mu}+\frac{\sigma(1-\rho)\alpha}{(\lambda_0+\alpha+\mu)(\lambda_0+\gamma_2+\varphi_2+\mu)}+\frac{\varepsilon\omega}{\lambda_0+\gamma_3+\mu}\right]\right|\\&\leq \frac{\beta S_0}{\sigma+\varepsilon+\mu}\left[\left|\frac{\sigma}{\lambda_0+\alpha+\mu}\right|+\left|\frac{\sigma(1-\rho)\alpha}{(\lambda_0+\alpha+\mu)(\lambda_0+\gamma_2+\varphi_2+\mu)}\right|+\left|\frac{\varepsilon\omega}{\lambda_0+\gamma_3+\mu}\right|\right]\\&\leq \frac{\beta S_0}{\sigma+\varepsilon+\mu}\left[\left|\frac{\sigma}{\alpha+\mu}\right|+\left|\frac{\sigma(1-\rho)\alpha}{(\alpha+\mu)(\gamma_2+\varphi_2+\mu)}\right|+\left|\frac{\varepsilon\omega}{\gamma_3+\mu}\right|\right] \quad (5)\\&= \frac{\beta S_0}{\sigma+\varepsilon+\mu}\left[\frac{\sigma}{\alpha+\mu}+\frac{\sigma(1-\rho)\alpha}{(\alpha+\mu)(\gamma_2+\varphi_2+\mu)}+\frac{\varepsilon\omega}{\gamma_3+\mu}\right]\\&= \mathcal{R}_C.\end{aligned}$$

Therefore, in the case of $\mathcal{R}_C < 1$, it can be known from inequality (5) that

$$\operatorname{Re}[F(\lambda_0)] < 1. \tag{6}$$

At this point, it can be seen that there is a contradiction between inequalities (4) and (6), so the initial assumption about $\lambda_0$ is invalid, that is, all the characteristic roots of the characteristic equation $\det(\lambda_0 E - M_1) = 0$ have negative real parts. That is to say, when $\mathcal{R}_C < 1$, the disease-free equilibrium point $P_0$ of system (1) is locally asymptotically stable.

Theorem 2.4.2. When $\mathcal{R}_C < 1$, the disease-free equilibrium point $P_0$ of system (1) is globally asymptotically stable.

Define function

$$h(x) = x - 1 - \ln x, x > 0.$$

For any $x > 0$ there is $h(x) \geq 0$, and $h(x) = 0$ if and only if $x = 1$.

Proof. Let $U(t) = (S(t), E_1(t), E_2(t), I_1(t), I_2(t), A(t), R(t))$ be the solution of system (1). Clearly $U(t)$ is non-negative and uniformly bounded, and $S(t) > 0$ holds for any $t \geq 0$.

Define a Lyapunov function

$$V(t) = S_0 h\left(\frac{S(t)}{S_0}\right) + \mathcal{R}_C E_1(t) + \frac{\beta S_0}{\alpha + \mu} E_2(t) + \frac{\omega \beta S_0}{\gamma_3 + \mu} A(t)$$

$$+ \frac{\beta S_0}{\gamma_2 + \varphi_2 + \mu}[(1-\rho)E_2(t) + I_2(t)] - \frac{\beta S_0 (1-\rho)\mu}{(\alpha + \mu)(\gamma_2 + \varphi_2 + \mu)} E_2(t).$$

The derivative of $V(t)$ with respect to $t (t > 0)$ along the system (1) is

$$\frac{dV(t)}{dt} = \left(1 - \frac{S_0}{S(t)}\right)\left[\Lambda - \beta S(t)(E_2(t) + I_2(t) + \omega A(t)) - \mu S(t)\right]$$

$$+ \mathcal{R}_C\left[\beta S(t)(E_2(t) + I_2(t) + \omega A(t)) - (\sigma + \varepsilon + \mu)E_1(t)\right]$$

$$+ \frac{\beta S_0}{\alpha + \mu}\left[\sigma E_1(t) - (\alpha + \mu)E_2(t)\right] + \frac{\omega \beta S_0}{\gamma_3 + \mu}\left[\varepsilon E_1(t) - (\gamma_3 + \mu)A(t)\right]$$

$$+ \frac{\beta S_0}{\gamma_2 + \varphi_2 + \mu}\left[\sigma(1-\rho)E_1(t) - \mu(1-\rho)E_2(t) - (\gamma_2 + \varphi_2 + \mu)I_2(t)\right]$$

$$- \frac{\beta S_0(1-\rho)\mu}{(\alpha + \mu)(\gamma_2 + \varphi_2 + \mu)}\left[\sigma E_1(t) - (\alpha + \mu)E_2(t)\right]$$

$$= -\frac{\mu}{S(t)}\left[S(t) - S_0\right]^2 + (\mathcal{R}_C - 1)\beta S(t)(E_2(t) + I_2(t) + \omega A(t)).$$

When $\mathcal{R}_C < 1$, there are $\frac{dV(t)}{dt} \leq 0$, and $\frac{dV(t)}{dt} = 0$ if and only if $U(t) = P_0$. Therefore, the largest invariant set of $\left\{U(t) \bigg| \frac{dV(t)}{dt} = 0\right\}$ is $\{P_0\}$. According to the LaSalle's Invariance Principle, when $\mathcal{R}_C < 1$, the disease-free equilibrium point $P_0$ of system (1) is globally asymptotically stable.

Theorem 2.4.3. When $\mathcal{R}_C > 1$, the disease-free equilibrium point $P_0$ of system (1) is unstable.

Proof. Similarly, when $\mathcal{R}_C > 1$, assume that there is $\lambda_1$ satisfying $\text{Re}(\lambda_1) \geq 0$, such that $\det(\lambda_1 E - M_1) = 0$. Expanding the characteristic equation, we can obtain

$$(\lambda_1 + \mu)^2(\lambda_1 + \alpha + \mu)(\lambda_1 + \gamma_1 + \varphi_1 + \mu)(\lambda_1 + \gamma_2 + \varphi_2 + \mu)(\lambda_1 + \gamma_3 + \mu)$$
$$\{\lambda_1 + (\sigma + \varepsilon + \mu)[1 - F(\lambda_1)]\} = 0, \qquad (7)$$

and

$$F(\lambda_1) = \frac{\beta S_0}{\sigma + \varepsilon + \mu}\left[\frac{\sigma}{\lambda_1 + \alpha + \mu} + \frac{\sigma(1-\rho)\alpha}{(\lambda_1 + \alpha + \mu)(\lambda_1 + \gamma_2 + \varphi_2 + \mu)} + \frac{\varepsilon \omega}{\lambda_1 + \gamma_3 + \mu}\right].$$

The characteristic equation (7) is equivalent to

$$(\lambda_1 + \alpha + \mu)(\lambda_1 + \gamma_2 + \varphi_2 + \mu)(\lambda_1 + \gamma_3 + \mu)\{\lambda_1 + (\sigma + \varepsilon + \mu)[1 - F(\lambda_1)]\} = 0. \qquad (8)$$

To simplify the subsequent formulation, redefine the notation as

$$B_1 = \alpha + \mu, B_2 = \gamma_2 + \varphi_2 + \mu, B_3 = \gamma_3 + \mu, C_1 = \sigma + \varepsilon + \mu,$$
$$C_2 = \sigma, C_3 = \sigma(1-\rho)\alpha, C_4 = \varepsilon\omega, D = \beta S_0.$$

Equation (8) is now rewritten as

$$(\lambda_1 + B_1)(\lambda_1 + B_2)(\lambda_1 + B_3)\left\{\lambda_1 + C_1\left\{1 - \frac{D}{C_1}\left[\frac{C_2}{\lambda_1 + B_1} + \frac{C_3}{(\lambda_1 + B_1)(\lambda_1 + B_2)} + \frac{C_4}{\lambda_1 + B_3}\right]\right\}\right\} = 0.$$

The result of expanding the above formula is

$$\lambda_1^4 + a_1\lambda_1^3 + a_2\lambda_1^2 + a_3\lambda_1 + a_4 = 0,$$

Here

$a_1 = B_1 + B_2 + B_3 + C_1,$
$a_2 = B_1B_2 + B_1B_3 + B_2B_3 + B_1C_1 + B_2C_1 + B_3C_1 - DC_2 - DC_4,$
$a_3 = B_1B_2B_3 + B_1B_2C_1 + B_1B_3C_1 + B_2B_3C_1 - DB_2C_2 - DB_3C_2 - DC_3 - DB_1C_4 - DB_2C_4,$
$a_4 = B_1B_2B_3C_1 - DB_2B_3C_2 - DB_1B_2C_4 - DB_3C_3.$

Define function

$$G(\lambda) = \lambda^4 + a_1\lambda^3 + a_2\lambda^2 + a_3\lambda + a_4.$$

Under the new notation, the number of controlled regenerations is

$$\mathcal{R}_C = \frac{D}{C_1}\left(\frac{C_2}{B_1} + \frac{C_3}{B_1B_2} + \frac{C_4}{B_3}\right).$$

Therefore, when $\mathcal{R}_C > 1$, there is

$$G(0) = B_1B_2B_3C_1 - DB_2B_3C_2 - DB_1B_2C_4 - DB_3C_3 = \frac{1}{B_1B_2B_3C_1}(1 - \mathcal{R}_C) < 0.$$

In addition, there is $G(+\infty) = +\infty$, so equation (8) has at least one positive root. That is to say, the characteristic equation $\det(\lambda_1 E - M_1) = 0$ has at least one positive root, that is, the assumption is true. This indicates that, when $\mathcal{R}_C > 1$, the disease-free equilibrium point $P_0$ of system (1) is unstable.

Theorem 2.4.4. When $\mathcal{R}_C > 1$, the system (1) is locally asymptotically stable at the endemic equilibrium point $P^*$.

Proof. Let the Jacobian matrix of system (1) at the endemic equilibrium point $P^*$ be $M_2$. When $\mathcal{R}_C > 1$, assume that there is $\lambda_2$ satisfying $\text{Re}(\lambda_2) \geq 0$, such that $\det(\lambda_2 E - M_2) = 0$. Expanding the characteristic equation, we can obtain

$$(\lambda_2 + \mu)(\lambda_2 + \alpha + \mu)(\lambda_2 + \gamma_1 + \varphi_1 + \mu)(\lambda_2 + \gamma_2 + \varphi_2 + \mu)(\lambda_2 + \gamma_3 + \mu)$$
$$\left\{(\lambda_2 + \mu)\left[\lambda_2 + (\sigma + \varepsilon + \mu)\left(1 - \frac{S^*}{S_0}F(\lambda_2)\right)\right] + \beta\left(E_2^* + I_2^* + \omega A^*\right)(\lambda_2 + \sigma + \varepsilon + \mu)\right\} = 0, \quad (9)$$

where,

$$F(\lambda_2) = \frac{\beta S_0}{\sigma + \varepsilon + \mu}\left[\frac{\sigma}{\lambda_2 + \alpha + \mu} + \frac{\sigma(1-\rho)\alpha}{(\lambda_2 + \alpha + \mu)(\lambda_2 + \gamma_2 + \varphi_2 + \mu)} + \frac{\varepsilon\omega}{\lambda_2 + \gamma_3 + \mu}\right].$$

The characteristic equation (9) is equivalent to

$$(\lambda_2 + \mu)(\lambda_2 + \sigma + \varepsilon + \mu) + \beta\left(E_2^* + I_2^* + \omega A^*\right)(\lambda_2 + \sigma + \varepsilon + \mu)$$
$$= \frac{S^*(\sigma + \varepsilon + \mu)}{S_0}(\lambda_2 + \mu)F(\lambda_2).$$

Using equation (2), the above equation can be transformed into

$$(\lambda_2 + \mu)(\lambda_2 + \sigma + \varepsilon + \mu) + \frac{\mathcal{R}_C E_1^*}{S_0}(\sigma + \varepsilon + \mu)(\lambda_2 + \sigma + \varepsilon + \mu)$$
$$= \frac{S^*(\sigma + \varepsilon + \mu)}{S_0}(\lambda_2 + \mu)F(\lambda_2).$$

The above formula can be obtained by simplification as

$$1 + \frac{\mathcal{R}_C E_1^*(\sigma + \varepsilon + \mu)}{S_0} \cdot \frac{1}{\lambda_2 + \mu} = \frac{S^*(\sigma + \varepsilon + \mu)}{S_0(\lambda_2 + \sigma + \varepsilon + \mu)}F(\lambda_2). \quad (10)$$

Furthermore, from the equilibrium equation we can obtain

$$\beta S^*\left(E_2^* + I_2^* + \omega A^*\right) = (\sigma + \varepsilon + \mu)E_1^*.$$

Using equation (2) to replace the parameters of the above equation, we can obtain

$$\beta S^*\left(\frac{\sigma}{\alpha + \mu} + \frac{\sigma(1-\rho)\alpha}{(\alpha + \mu)(\gamma_2 + \varphi_2 + \mu)} + \frac{\omega\varepsilon}{\gamma_3 + \mu}\right)E_1^* = (\sigma + \varepsilon + \mu)E_1^*,$$

and finally, we obtain

$$\mathcal{R}_C = \frac{S_0}{S^*}.$$

Similar to the proof of inequality (5), taking the norm on the right side of the equal sign of equation (10) can get

$$\left|\frac{S^*(\sigma + \varepsilon + \mu)}{S_0(\lambda_2 + \sigma + \varepsilon + \mu)}F(\lambda_2)\right| \leq \left|\frac{S^*(\sigma + \varepsilon + \mu)}{S_0(\lambda_2 + \sigma + \varepsilon + \mu)}\mathcal{R}_C\right| = \left|\frac{\sigma + \varepsilon + \mu}{\lambda_2 + \sigma + \varepsilon + \mu}\right| \leq 1. \quad (11)$$

Since we assume $\text{Re}(\lambda_2) \geq 0$, there is $\text{Re}[1/(\lambda_2 + \mu)] > 0$, and the left-hand side of equation (10) satisfies

$$\text{Re}\left(1 + \frac{\mathcal{R}_C E_1^*(\sigma + \varepsilon + \mu)}{S_0} \cdot \frac{1}{\lambda_2 + \mu}\right) = 1 + \frac{\mathcal{R}_C E_1^*(\sigma + \varepsilon + \mu)}{S_0}\text{Re}\left(\frac{1}{\lambda_2 + \mu}\right) > 1.$$

Therefore, taking the norm on the left side of the equal sign in equation (10), we get

$$\left| 1 + \frac{\mathcal{R}_C E_1^* (\sigma + \varepsilon + \mu)}{S_0} \cdot \frac{1}{\lambda_2 + \mu} \right| > 1. \tag{12}$$

At this point we can see that the inequalities (11) and (12) are contradictory to the equation (10), so the initial assumption about $\lambda_2$ is invalid, that is, all the eigenvalues of the characteristic equation $\det(\lambda_2 E - M_2) = 0$ have negative real parts. That is, when $\mathcal{R}_C > 1$, the system (1) is locally asymptotically stable at the endemic equilibrium point $P^*$.

3. Numerical fitting

3.1 Selection of fitting objects

The SARS-CoV-2 has mutated several times since it was first discovered in late 2019. Among them, there are 7 variants of concern that caused a global pandemic, including 614D, 614G, Alpha, Beta, Gamma, Delta and Omicron variants [11]. 614D is an original variant and there is a lack of data on its infectors; although the Beta variant has strong immune evasion ability and its transmission is stronger than that of the Alpha variant, its epidemic process is quickly overtaken by the more transmissible variant; similarly, the Gamma variant, although more pathogenic, was rapidly replaced by the Delta variant in its epidemic process due to its weaker transmission; therefore, we do not consider the above three variants. This study focuses only on four of the major virus variants that are widely circulating worldwide, namely 614G, Alpha, Delta and Omicron variants.

For the 614G variant, we select the epidemic data from June 1, 2020 to July 31, 2020 in England for simulation, which was recorded as the first time period. The COVID-19 outbreak in England has since progressed slowly until early November 2020, when a second nationwide lockdown began and lasted for about a month. The second national lockdown did play a relatively obvious role in curbing the spread of the epidemic, but due to the dual effects of the gradual relaxation of the lockdown policy and the emergence of the Alpha variant, England has experienced the third large-scale outbreak in a short period of time [12]. Therefore, for the Alpha variant, we select the epidemic data from December 1, 2020 to December 31, 2020 in England for simulation, which was recorded as the second time period.

Since then, England has entered the third national lockdown and lasted for about two months, which has brought the epidemic in England under certain control. However, similarly, England began its third large-scale outbreak in May 2021 due to the dual effects of the gradual easing of the third lockdown policy and the emergence of the Delta variant. Since the Delta variant first appeared in England in mid-April 2021, it has replaced the Alpha variant as the predominant variant in England in just over a month. Therefore, for the Delta variant, we select the epidemic data from June 1, 2021 to June 31, 2021 in England for simulation, which is recorded as the third time period. Since then, the UK has been adopting a lax epidemic prevention policy, which has kept the daily new case data at a high position.

For the Omicron variant, we select the epidemic data from March 18, 2022 to April 9, 2022 in Shanghai, China for simulation, which is recorded as the fourth time period [13]. There are two main reasons for choosing the epidemic data in Shanghai during this period as the research object. First, this round of epidemic in Shanghai is basically caused by a single Omicron variant, which facilitates our comparative analysis. Second, this time period this is the early stage of the outbreak in Shanghai, will be used to fit the mathematical model. The relative time stage of the fitting time period selected relative to this round of epidemics are shown in Figure 2. The first three fitting time periods we selected are the early stages of three large-scale outbreaks in England, and the fourth time period is the early stages of the Shanghai outbreak in April 2022. We use the established model to fit the epidemic data in these four time periods, and comparatively analyze the role of asymptomatic infections in the spread of the epidemic.

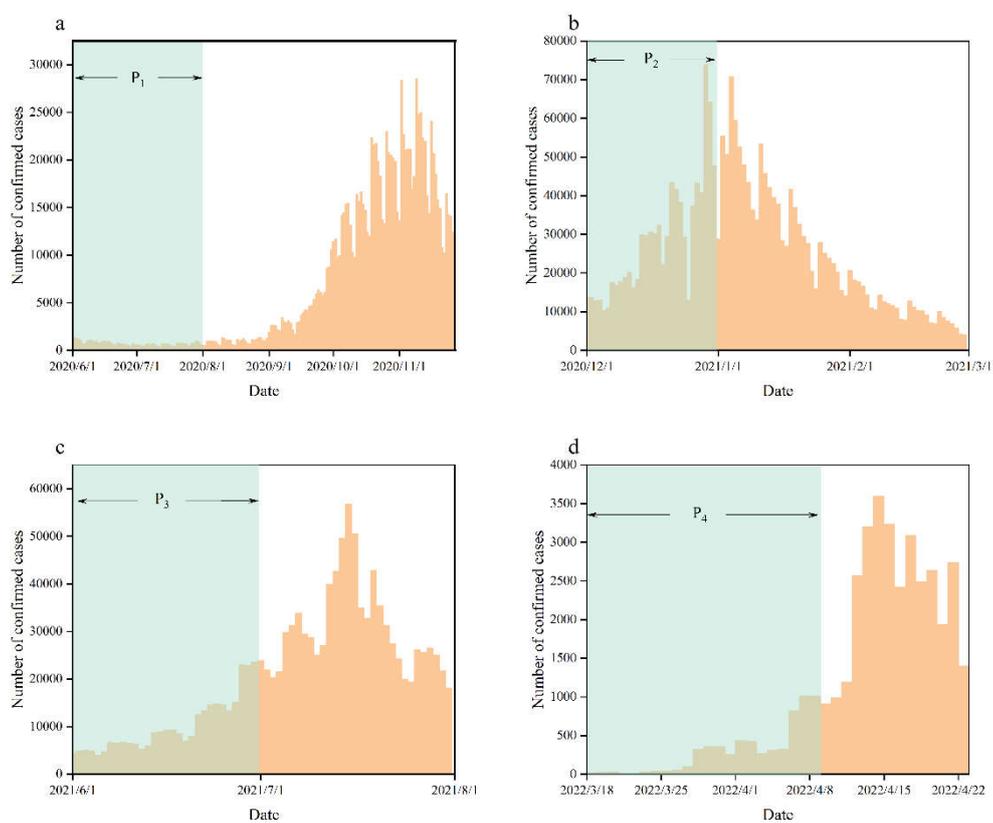

**Fig.2 Relative position map of the four data fitting time periods during the corresponding outbreak. a**, **b** and **c**, The relative position map of the selected first, second and third data fitting time periods during the first three outbreaks of the corresponding epidemic in England. **d**, The relative position map of the selected fourth data fitting time period during the corresponding outbreak in Shanghai in April 2022.

3.2 Numerical fitting results

For the selected 614G, Alpha, Delta and Omicron variants, we simulate the epidemic data of different periods using the constructed mathematical model. Among them, the epidemic data corresponding to the 614G, Alpha and Delta variants are all from the official website of the British government (https://coronavirus.data.gov.uk/), and the epidemic data corresponding to the Omicron variant are from the Shanghai Municipal Health Commission (https://wsjkw.sh.gov.cn/).

The fitting results of the epidemic data in the four time periods are shown in Figure 3, in which Figures 3a to 3d show the data fitting in the first to fourth time periods in turn. From the figure, we can see that the simulation results of the epidemic data in the four time periods are relatively satisfactory, which indicates that our model has good universality. Due to the small value and large fluctuation of the original data, the simulation results of the first time period and the fourth time period are not perfect in some parts, but they are generally acceptable. The related dynamic parameters obtained by fitting are shown in Table 1.

Table.1 Fitting results of dynamic parameters in model (1).

| Parameter | Value | | | | Reference |
|---|---|---|---|---|---|
| | 614G | Alpha | Delta | Omicron | |
| $\Lambda$ | 1740 | 1740 | 1740 | 216 | [14,15] |
| $\mu$ | $2.5753 \times 10^{-5}$ | $2.5753 \times 10^{-5}$ | $2.5753 \times 10^{-5}$ | $2.4303 \times 10^{-5}$ | [14,15] |
| $\beta$ | $5.3720 \times 10^{-9}$ | $7.2151 \times 10^{-9}$ | $9.0205 \times 10^{-9}$ | $3.2493 \times 10^{-8}$ | Fitted |
| $\sigma$ | 0.1975 | 0.1975 | 0.1975 | 0.1975 | [16] |
| $\varepsilon$ | 0.3415 | 0.5748 | 0.6768 | 0.5745 | Fitted |
| $\alpha$ | 0.5000 | 0.5000 | 0.5000 | 0.5000 | [17] |
| $\omega$ | 0.6524 | 0.6524 | 0.6524 | 0.6524 | [7] |
| $\rho$ | 0.4689 | 0.1103 | 0.2266 | 0.5266 | Fitted |
| $\gamma_1$ | 0.0588 | 0.0588 | 0.0588 | 0.0588 | [18-20] |
| $\gamma_2$ | 0.0769 | 0.0811 | 0.0704 | 0.0537 | Fitted |
| $\gamma_3$ | 0.2770 | 0.3746 | 0.4810 | 0.4149 | Fitted |
| $\varphi_1$ | $1.7826 \times 10^{-5}$ | $1.7826 \times 10^{-5}$ | $1.7826 \times 10^{-5}$ | $1.7826 \times 10^{-5}$ | [21] |
| $\varphi_2$ | $5.5963 \times 10^{-3}$ | $4.4054 \times 10^{-3}$ | $5.0410 \times 10^{-3}$ | $5.0179 \times 10^{-3}$ | Fitted |

3.3 Horizontal comparison of fitting results

From the perspective of mathematical model simulation, we verify the characteristics of enhanced infectivity and weakened symptoms of infected persons of the four SARS-CoV-2 variants concerned in this study, as shown in Figure 4. We count the proportions of the fitting results of the three types of infected persons. From the circular distribution map of the types of infected people, we can see that the proportion of asymptomatic patients in the four time periods was 38.05%, 52.89%, 54.02% and 62.03% respectively, showing a clear upward trend. This result shows that the overall symptom level of the infected person is decreasing, which is consistent with the clinical statistics of the decrease in the severe rate of the infected person.

In addition, we also compare the propagation speed and propagation strength of different variants. From the perspective of the model, we predict the subsequent development of the epidemic in four time periods. From Figures 4c and 4f, we can see that the peak time of the epidemic development process caused by the variants of interest in this study gradually shortens, and the peak value of daily confirmed cases gradually increases, and there is even an order of magnitude difference. This indicates that the speed of transmission, as well as the intensity of transmission, of the four COVID-19 variants of interest in this study has increased. In addition, it can also be seen from the dynamic parameters obtained by fitting that the basic propagation rate has increased significantly.

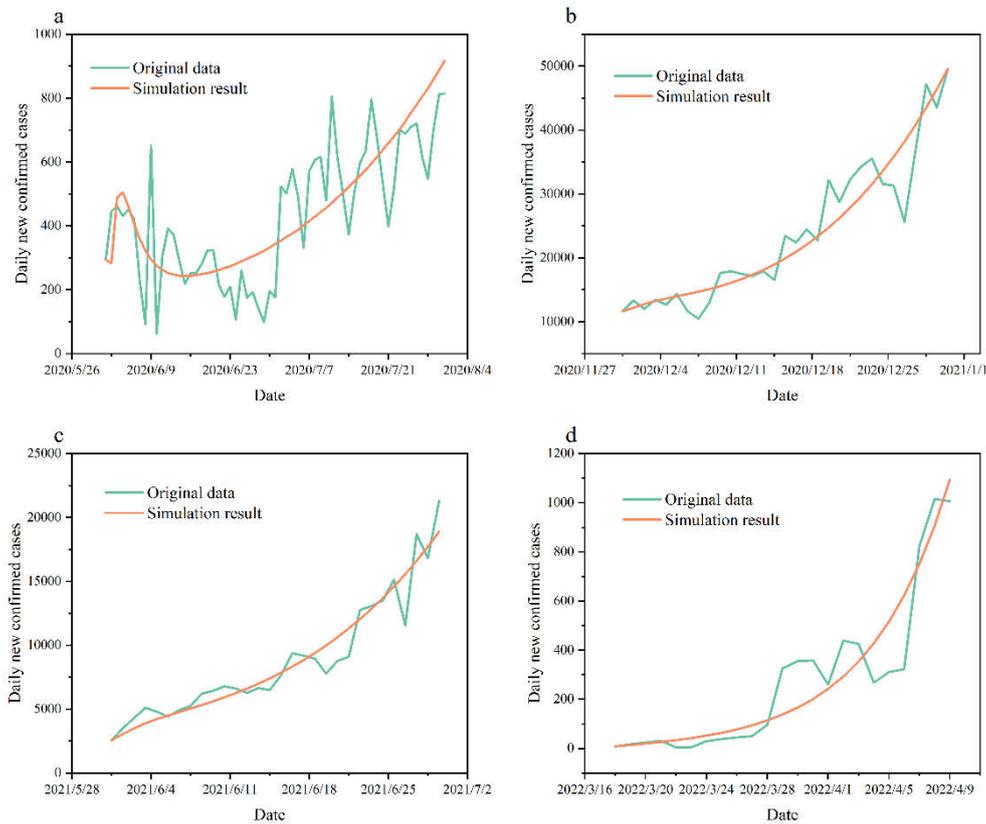

**Fig.3 Simulation results of the daily number of newly confirmed cases of epidemic data in four time periods. a**, Simulation results in the first time period, corresponding to the 614G variant. **b**, Simulation results in the second time period, corresponding to the Alpha variant. **c**, Simulation results in the third time period, corresponding to the Delta variant. **d**, Simulation results in the fourth time period, corresponding to the Omicron variant.

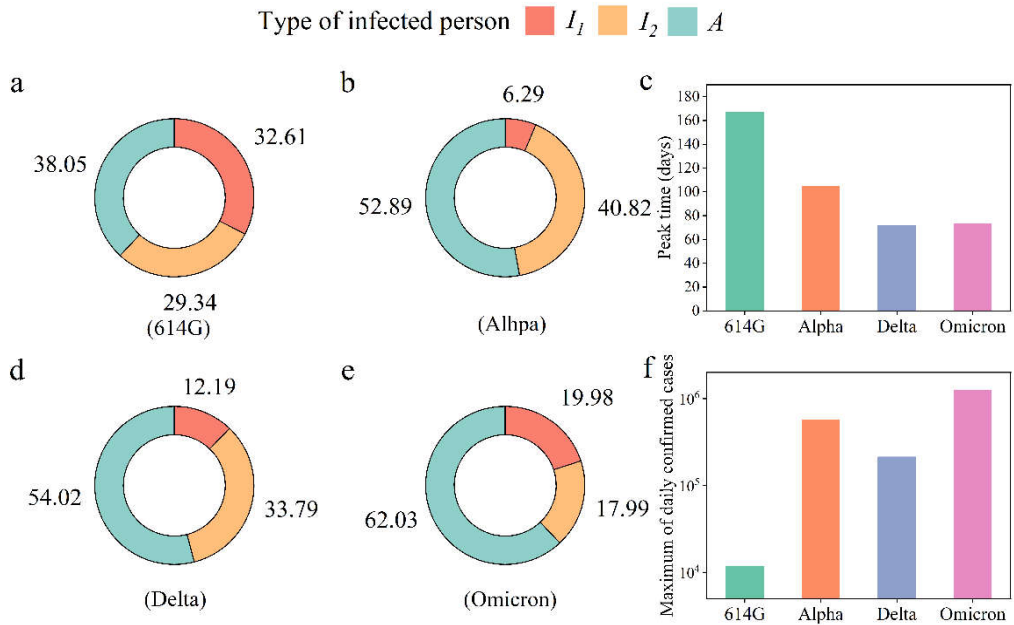

**Fig.4 Horizontal comparison results of epidemic data simulation in four time periods. a**, **b**, **d** and **e**, Changes in the proportion of the three types of infected persons. **c**, Model prediction results of the peak time of the epidemic in four time periods. **f**, Model predictions of maximum daily confirmed cases over four time periods.

3.4 The impact of the rate of participation in testing among symptomatic infected persons on the spread of the epidemic

Although the rate of testing among those with symptomatic infection is fixed in the model simulation, in practice, with the increase in the number of newly diagnosed patients every day, certain prevention and control measures will inevitably be taken to curb the spread of the epidemic. The most direct and effective prevention and control measure is to increase the detection intensity, which is reflected in the model as the increase of the parameter $\rho$. Therefore, we explore the effect of the change of parameter $\rho$ on the spread of the epidemic.

Based on the simulation results over four time periods, while keeping other parameters unchanged, re-prediction is carried out for the cases of $\rho = 20\%, 40\%, 60\%$, and $80\%$, respectively, as shown in Figure 5 and Figure 6. From the results, we can see that with the increase of $\rho$, the cumulative number of infected people has dropped significantly. In addition, the proportion of asymptomatic infections has declined to a certain extent, but the decline is not significant. As the number of symptomatic infected people tested increases, the impact of symptomatic infection on the spread of the epidemic gradually decreases. Since the model assumes that asymptomatic infected people will not participate in the detection, although the transmission intensity is lower than that of symptomatic infected people, it will still cause a certain amount of transmission. Therefore, even if the proportion of symptomatic infected people participating in the test reaches 80%, there are still a certain proportion of asymptomatic infections.

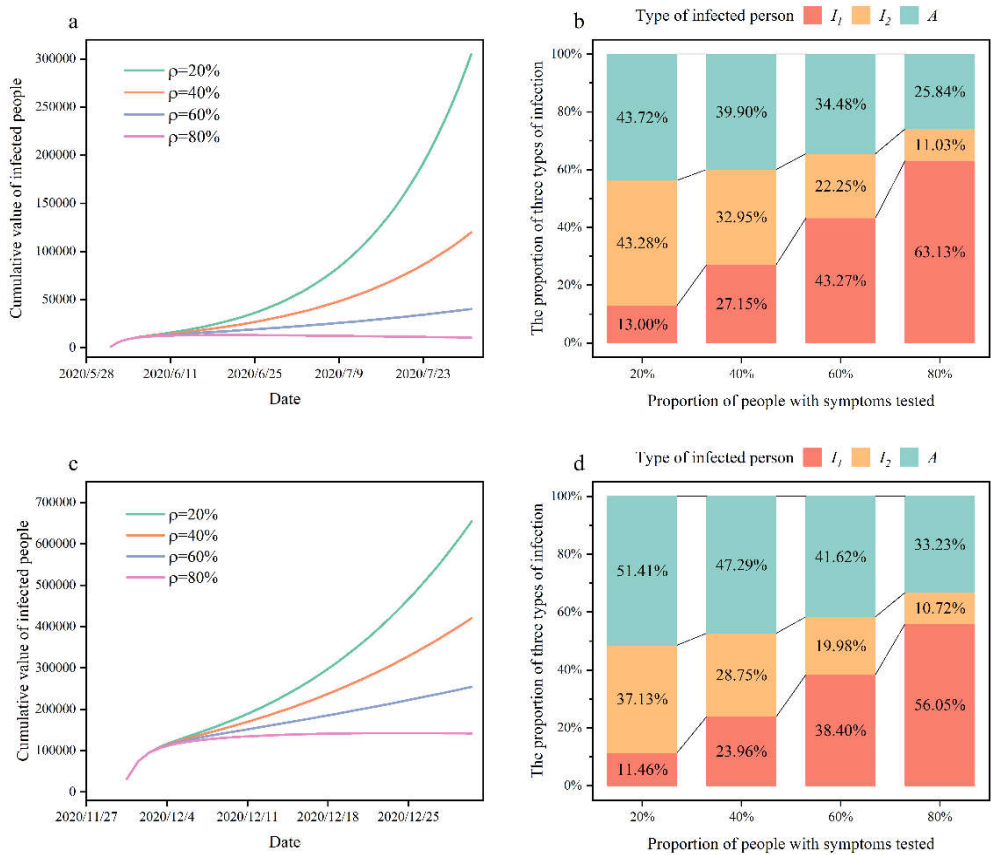

**Fig.5 The impact of the tested ratio ρ among symptomatic infected persons on the spread of the epidemic. a**, The influence of parameter ρ changes on the cumulative amount of infected people, based on the fitting results of the first time period. **b**, The impact of parameter ρ changes on the proportion of different types of infected people, based on the fitting results of the first time period. **c**, The influence of parameter ρ changes on the cumulative amount of infected people, based on the fitting results of the second time period. **d**, The impact of parameter ρ changes on the proportion of different types of infected people, based on the fitting results of the second time period.

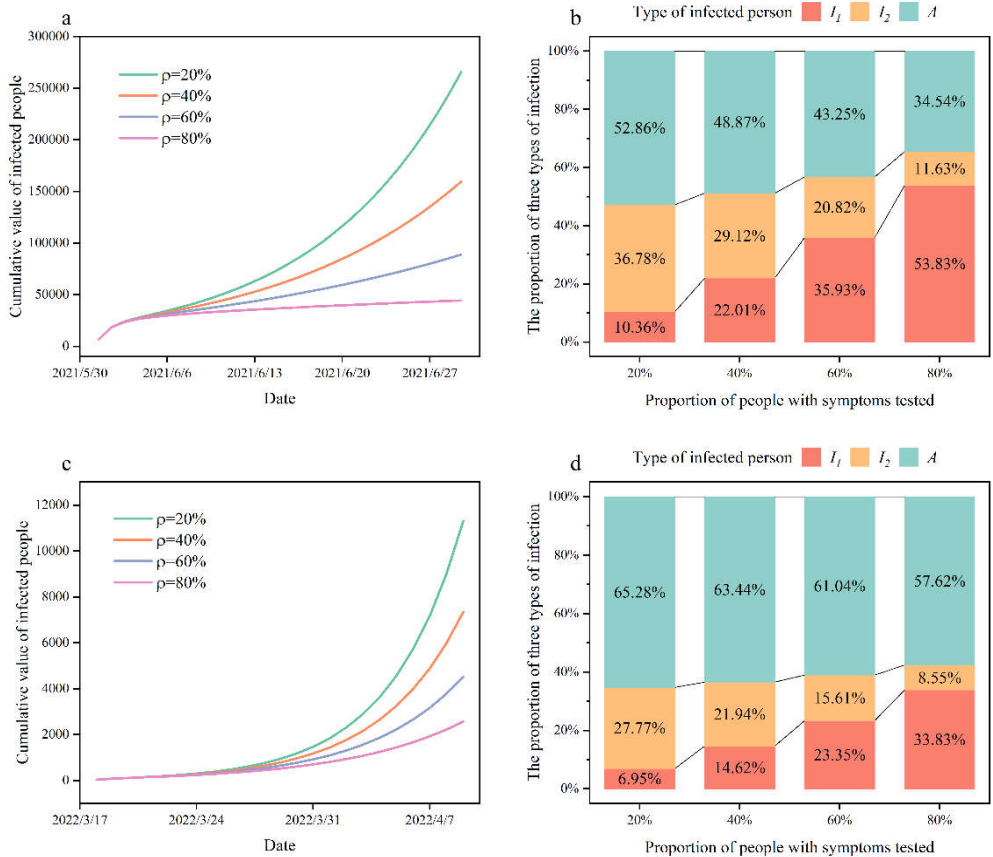

**Fig.6 The impact of the tested ratio ρ among symptomatic infected persons on the spread of the epidemic. a**, The influence of parameter ρ changes on the cumulative amount of infected people, based on the fitting results of the third time period. **b**, The impact of parameter ρ changes on the proportion of different types of infected people, based on the fitting results of the third time period. **c**, The influence of parameter ρ changes on the cumulative amount of infected people, based on the fitting results of the fourth time period. **d**, The impact of parameter ρ changes on the proportion of different types of infected people, based on the fitting results of the fourth time period.

In addition, we compare the differences between the different variants when the rate of detection in symptomatic infected individuals changed. As can be seen from Figure 7, when the detection ratio of symptomatic infections increases from 20% to 80%, the decline ratio of the total number of infected people decreases slowly with the variation of variants. However, the decline rates remain above 70%, indicating that increased testing efforts can still effectively reduce the total number of infected people. In addition, the proportion of asymptomatic infections decrease significantly with the change of variants. From Figure 7, we can see that for the Omicron variant, when the detection ratio increased from 20% to 80%, the proportion of asymptomatic infections drop by only about 10%. The specific data are shown in Table 2. From the results of this part, we can see that in the process of epidemic prevention, only improving the detection of symptomatic infections cannot completely block the spread of the epidemic. It is necessary to cooperate with certain social isolation policies to effectively block hidden transmission caused by asymptomatic infections.

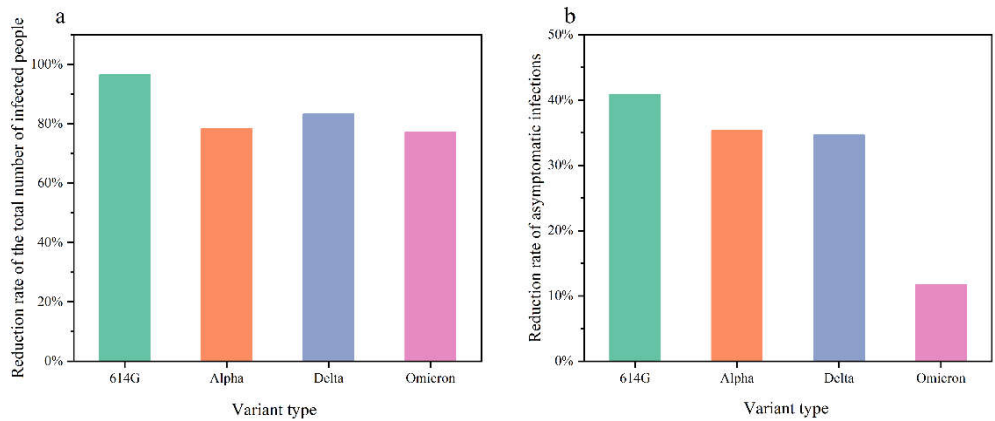

**Fig.7 The impact of the tested ratio ρ among symptomatic infected persons on the spread of the epidemic. a**, The reduction rate in the total number of infected people when the detection rate goes from 20% to 80%. **b**, The reduction rate in asymptomatic infections when testing rates go from 20% to 80%.

Table.2 The impact of parameter ρ changes on the spread of the epidemic in different time periods.

|  | 614G | Alpha | Delta | Omicron |
|---|---|---|---|---|
| The percentage of decline in the total number of infected people | 96.51% | 78.37% | 83.31% | 77.30% |
| The percentage of decline in the number of asymptomatic infections | 40.90% | 35.37% | 34.66% | 11.72% |

3.5 Model-based prediction results

We use the established model to predict the spread of the epidemic. Since the English government no longer publishes epidemic data from May 2022, we choose WHO Coronavirus (COVID-19) Dashboard (https://covid19.who.int/) as the source of epidemic data to fit and predict the development of the epidemic in the UK. As shown in Figure 8, we forecast the development of the epidemic in the UK from late 2022 to early 2023. Predictions suggest a larger outbreak in the UK in the winter of 2022, which would last for around four months. In order to avoid the collapse of the healthcare system and the negative impact on productive life, the UK should take the necessary measures to hinder the spread of the epidemic.

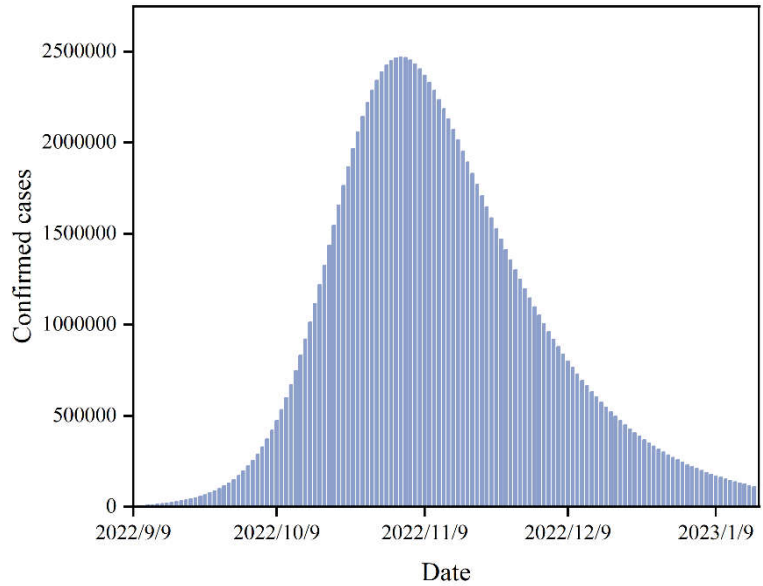

**Fig.8 The latest UK epidemic fitting and forecast results.** Based on the predicted results, the outbreak will last from September 2022 to January 2023. The number of confirmed cases will peak in early November 2022.

4. Conclusion

In this study, a SEIAR infectious disease dynamics model was established. We first analyze the local stability at the equilibrium point of the model system. As a result, when $\mathcal{R}_C < 1$, the disease-free equilibrium point $P_0$ of the system is locally asymptotically stable; when $\mathcal{R}_C > 1$, the disease-free equilibrium point $P_0$ of the system is unstable and the endemic equilibrium point $P^*$ is locally asymptotically stable. In addition, we use the method of constructing the Lyapunov function to prove that the disease-free equilibrium point $P_0$ of the system is globally asymptotically stable when $\mathcal{R}_C < 1$. In addition, we use the model of the architecture to numerically fit the COVID-19 epidemic data. We select corresponding data for the four main variants, 614G, Alpha, Delta and Omicron, and performed fitting analysis. The comparative analysis of the fitting results shows that the proportion of asymptomatic infections among infected persons is gradually increasing, and the propagation speed and propagation strength of the variant are increased to a certain extent. Finally, we analyze the impact of the detection ratio of symptomatic infections on the spread of the epidemic. With the increase of the detection ratio, the cumulative cases drop significantly, but asymptomatic infections still maintain a certain proportion.

## 5. Discussion

Since asymptomatic infections are difficult to be actively detected, the transmission caused by asymptomatic infections is difficult to be interrupted under relatively loose epidemic prevention policies. Especially after the end of the last round of the epidemic, the epidemic prevention policy tends to be gradually loosened as the number of daily confirmed cases declines, which will cause the number of asymptomatic infections to accumulate again. It can also be seen from the comparison results of this study that in general, the pathogenicity of the SARS-CoV-2 variants is decreasing but the infectivity is increasing, which means that the role of asymptomatic infected persons in the spread of the epidemic has increased significantly. And due to increased infectivity, the absolute number of symptomatic infections will not necessarily decrease. This shows that epidemic prevention and control measures are still very important, otherwise there will be continuous challenges to the medical carrying capacity.

Large-scale social activities will also cause large-scale spread of the epidemic in the short term. Taking the Tokyo Olympics as an example, during the entire event period, not only the cities where the competition were held, but also other cities that did not hold the competition experienced outbreaks [22]. During the Tokyo Olympics, Japan's epidemic prevention policy was relaxed, and social mobility increased. This made a large number of asymptomatic infections accumulated in the crowd, resulting in a full-scale outbreak of the epidemic. Therefore, when holding large-scale social activities, it is more necessary to strengthen the control of social mobility to reduce the probability of large-scale outbreaks in the short term.

In view of the increasing role of asymptomatic infections in the spread of the epidemic, non-pharmaceutical intervention (NPI) may play a more important role [23]. Common NPI measures include wearing masks, maintaining social distance, etc. Compared with symptomatic infected individuals, asymptomatic infected individuals have lower viral load and weaker transmission ability, which enables NPI to better block the transmission caused by asymptomatic infected individuals. In addition, vaccination policies still need to be actively implemented, as vaccination is the most effective way to reduce the risk of transmission overall. Therefore, in the process of epidemic prevention and control, when the detection capacity is limited, NPI and vaccination are important measures and means to avoid the collapse of the medical system.

**Acknowledgments**

This study was funded by Natural Science Foundation of China (NSFC 11871093, 11901027), Postgraduate Teaching Research and Quality Improvement Project of BUCEA (J2021010), BUCEA Post Graduate Innovation Project (PG2022139). We thank Professors Jingan Cui and Zhihang Peng for giving us much guidance, and all those who generously shared their time and materials for this study.